\documentclass[12pt]{article}
\usepackage{amssymb,latexsym}
\def\C{\mathbb C}
\def\R{\mathbb R}
\def\N{\mathbb N}
\def\const{{\mathrm{const}}}
\newtheorem{thm}{Theorem}[section]
\newtheorem{cor}{Corollary}[section]
\newtheorem{lem}{Lemma}[section]

\title{Real entire functions with real zeros and
a conjecture of Wiman}
\author{Walter Bergweiler\thanks{Supported by the German--Israeli
Foundation for Scientific Research and Development,
grant GIF G-643-117.6/1999, and by INTAS-99-00089.},
A. Eremenko\thanks{Supported by
the NSF grant DMS 0100512 and by the Humboldt
Foundation.}$\;$
and J.K. Langley\thanks{Supported by DAAD}}
\begin{document}
\maketitle
\section{Introduction}

An entire function is called real if it maps the real line into itself.
The main result of this paper is
\begin{thm}\label{thm0}
For every real entire function of infinite order with 
only real zeros, the second derivative
has infinitely many non-real zeros.
\end{thm}
This conclusion is not true for the first derivative as
the example $\exp(\sin z)$ shows. For real entire functions
with finitely many zeros, all of them real, Theorem \ref{thm0}
was proved in \cite{BF}.
Theorem \ref{thm0} can be considered as an extension to
functions of infinite order of the following
result of Sheil-Small
\cite{SS}, conjectured by
Wiman in 1914 \cite{Al1,Al2}. For every integer
$p\geq 0$, denote by $V_{2p}$ the set
of entire functions of the form
$$f(z)=\exp(-az^{2p+2})g(z),$$
where $a\geq 0$ and $g$ is a real entire function
with only 
real zeros of genus at most $2p+1$, and set
$U_0=V_0$ and
$U_{2p}=
V_{2p}\backslash V_{2p-2}$ for $p\geq 1$.
Thus the class of all real entire functions of finite order
with real zeros is represented as a union of disjoint
subclasses $U_{2p},\; p=0,1,\ldots.$
\vspace{.1in}

\noindent
{\bf Theorem A} (Sheil-Small) {\em If $f\in U_{2p}$ then
$f''$ has at least $2p$ non-real zeros.}
\vspace{.1in}

Applying Theorem \ref{thm0} to functions of the form
$$f(z)=\exp\int_0^zg(\zeta)\, d\zeta$$
we obtain
\begin{cor}\label{cor1}
For every real transcendental entire function $g$,
the function $g'+g^2$ has infinitely many non-real
zeros.
\end{cor}
For polynomials $g$ the corresponding result
was conjectured in \cite[Probl. 2.64 and 4.28]{Res}
and proved in
\cite{SS}: {\em If $g$ is a real polynomial
then $g'+g^2$ has at least $\deg g- 1$ non-real
zeros}. Corollary \ref{cor1} also follows from
the result of Bergweiler and Fuchs \cite{BF}. 

Theorems 1 and A together imply the following
\begin{cor}\label{cor2}
If $f$ is a real entire function
and $ff''$ has only real zeros then $f\in U_0$.
\end{cor}

We recall that $U_0$, the {\em Laguerre--P\'olya class}, coincides
with the closure of the set of all real
polynomials with only
real zeros, with respect to uniform convergence on
compact subsets of the plane.
This was proved by Laguerre \cite{Lag} for the case
of polynomials with positive zeros and 
by P\'olya \cite{Po2} in the general case.
It follows that $U_0$ is closed
under differentiation, so that all derivatives of a function
$f\in U_0$ have only real zeros. P\'olya \cite{Po2}
asked if the converse is true: {\em if all
derivatives of
a real entire function $f$ have only real zeros then
$f\in U_0$}.
This conjecture was proved by Hellerstein
and Williamson \cite{HelW1,HelW2}. More precisely, they showed
that {\em for a real entire function $f$, the condition
that $ff'f''$ has only real zeros implies $f\in U_0$.} 
Our Corollary \ref{cor2} shows that in this result
one can drop the assumption on the zeros of $f'$,
as Hellerstein and Williamson conjectured \cite[Probl. 2.64]{Res}. 

For the early history of results on the conjectures of
Wiman and P\'olya we refer to \cite{HelW1,LeO}, which contain
ample bibliography.
The main result of Levin and Ostrovskii \cite{LeO} is
\vspace{.1in}

\noindent
{\bf Theorem B} {\em If $f$ is a real entire
function
and all zeros of 
$ff''$ are real
then
\begin{equation}\label{00}
\log^+\log^+|f(z)|= O(|z|\log|z|),\quad z\to\infty.
\end{equation}}
This shows that a function satisfying the assumptions
of Corollary \ref{cor2} cannot grow too fast, but
there is a gap between Theorem B and Theorem A.
Our Theorem \ref{thm0} bridges this gap.

One important tool brought by Levin and Ostrovskii to
the subject was a factorization of the logarithmic
derivative of a real entire function $f$ with
only real zeros:
$$\frac{f'}{f}=\psi\phi,$$
where $\phi$ is a real entire function,
and either $\psi$ is a meromorphic function which maps
the upper half-plane $H=\{ z:\mathrm{Im}\, z>0\}$
into itself or $\psi\equiv 1$.
This factorization
was used in all subsequent work
in the subject.
A standard estimate for analytic functions mapping
the upper half-plane into itself shows that $\psi$
is neither too large nor too small away from the
real axis, so the asymptotic behavior of $f'/f$
mostly depends on that of $\phi$.
One can show that
$f$ is of finite order if and only if $\phi$ is
a polynomial.

The second major contribution of Levin and Ostrovskii was
the application of ideas from the value distribution 
theory of meromorphic functions \cite{GO,Haym1,Ne}.
Using Nevanlinna theory, Hayman \cite{Haym}
proved that for an entire function $f$, the condition
$f(z)f''(z)\neq 0,\; z\in\C$, implies that $f'/f$
is constant.
The assumptions of Theorem B mean that $f(z)f''(z)\neq 0$
in $H$. Levin and Ostrovskii adapted
Hayman's argument to functions in a half-plane to
produce an estimate for the logarithmic derivative.
An integration of this estimate gives (\ref{00}).
To estimate the logarithmic derivative using Hayman's
argument they applied an analogue of the
Nevanlinna characteristic
for meromorphic functions
in a half-plane, and proved an analogue of
the main technical result of Nevanlinna theory, the
lemma on the logarithmic derivative. This characteristic has two
independent origins, \cite{Le0} and \cite{Tsuji0},
and the name ``Tsuji characteristic''
was introduced in \cite{LeO}.

In this paper we use both main ingredients of the
work of Levin and Ostrovskii, the factorization of $f'/f$ 
and the Tsuji characteristic. 

Another important tool comes from  Sheil-Small's
proof of Theorem A.
His key idea was the study of topological
properties of the
auxiliary function
$$F(z)=z-\frac{f(z)}{f'(z)}.$$
In the last section
of his paper, Sheil-Small discusses the possibility of
extension of his method to functions of infinite order,
and proves the fact which turns out to be crucial:
{\em if $f$ is a real entire function, $ff''$ has only
real zeros, and $f'$ has a non-real zero,
then $F$ has a non-real asymptotic value.} 
In \S\ref{lem1pf} we prove
a generalization of this fact needed
in our argument. 

The auxiliary function $F$ appears when one solves the
equation $f(z)=0$ by Newton's method. This suggests
the idea of iterating $F$ and using the Fatou--Julia theory
of iteration of meromorphic functions. This was explored
by Eremenko and Hinkkanen, see, for example, \cite{Hinkk}.

Theorem \ref{thm0} will be proved by establishing 
a more general result
conjectured by Sheil-Small \cite{SS}.
Let $L$ be a real meromorphic function in the plane
with only simple poles, all of them real
and with positive residues.
It is known \cite{HelW1,LeO,SS}
that every such $L$ has a Levin--Ostrovskii representation
\begin{equation}
L=\psi\phi 
\label{1}
\end{equation}
in which:
\vspace{.1in}

(a) $\psi $ is meromorphic in the plane and real on the real axis;

(b) $\psi$ maps the upper
half-plane 
into itself, or
$\psi\equiv 1$;

(c) $\psi$ has a simple pole at every pole of $L$, and no other
poles;

(d) $\phi $ is a real entire function.
\vspace{.1in}

\noindent
We outline briefly how such a factorization
(\ref{1}) is obtained.
Let 
$$\ldots<a_{k-1}<a_k<a_{k+1}<\ldots$$
be the sequence of poles of $L$ enumerated in
increasing order.
The assumption that all poles are simple
and have positive residues implies
that there is at least one zero of $L$ in each
interval $(a_k,a_{k+1})$.
We choose one such zero in each interval and denote
it by $b_k$.
Then we set
$$\psi(z)=\prod_k\frac{1-z/b_k}{1-z/a_k},$$
with slight modifications if
$a_k b_k \leq 0$ for some $k$ or the set $\{ a_k \}$ is bounded above.
We then define $\phi$ by (\ref{1}), and
properties (a)--(d) follow (for the details see
\cite{HelW1,LeO,SS}). 


\begin{thm}\label{thm1}
Let $L$ be a function
meromorphic in the plane, real on the real axis, such that
all poles of $L$ are real, simple and have positive residues. Let
$\psi, \phi $ be as in
$(\ref{1})$ and {\rm (a), (b), (c), (d)}.
If $\phi$ is transcendental
then $L + L'/L$ has infinitely many non-real zeros.
\end{thm}

To deduce Theorem \ref{thm0} from
Theorem \ref{thm1} it suffices to note that if $f$ is a real entire function
with only real zeros then $L = f'/f$ is real meromorphic with only real
simple poles
and positive residues,
and thus has a representation (\ref{1}). Further,
$L + L'/L = f''/f'$. By an
argument of Hellerstein and Williamson
\cite[pp. 500-501]{HelW2}, the
function $\phi$ is transcendental if and only if $f$ has infinite
order.

\section{Preliminaries}

We will require the following well known consequence of 
Carleman's estimate 
for harmonic measure.

\begin{lem}\label{sub}
Let $u$ be a non-constant continuous subharmonic function in the plane.
For $r > 0$ let
$B(r, u) = \max \{ u(z) : |z| = r \}$, and let
$\theta (r)$ be the angular measure of that subset of the circle
$C(0, r) = \{ z \in \C : |z| = r \}$ on which $u(z) > 0$. Define
$\theta ^*(r)$ by $\theta ^*(r) =  \theta (r)$, except that
$\theta ^*(r) =  \infty $ if $u(z) > 0$ on the whole circle $C(0, r)$.
Then if $r > 2r_0$ and $B(r_0, u) > 1$ we have
$$
\log\| u^+(4re^{i\theta})\|\geq \log B(2r, u)-c_1 \geq
\int_{2r_0}^r \frac{\pi dt }{t \theta^*(t) } - c_2,
$$
in which $c_1$ and $c_2$ are absolute constants, and
$$\| u^+(re^{i\theta})\|=\frac{1}{2\pi}\int_{-\pi}^\pi
\max\{ u(re^{i\theta}),\,0\}\,d\theta.$$
\end{lem}

The first inequality follows from Poisson's formula,
and for the second we refer to \cite[Thm III.68]{Tsuji}. 
Note that in the case that $u=\log |f|$ where $f$ is an entire
function, $\| u^+(re^{i\theta})\|$ coincides with the
Nevanlinna characteristic $T(r,f)$.

Next, we need the characteristic function in a half-plane as developed
by Tsuji \cite{Tsuji0} and Levin and Ostrovskii \cite{LeO} (see also \cite{GO}
for a comprehensive treatment).
Let $f$ be a meromorphic function
in a domain containing the closed upper half-plane $\overline{H} =
\{ z \in \C : {\rm Im} (z) \geq  0 \}$.
For $t \geq 1$ let
$\mathfrak{n} (t, f) $ be the number of poles of $f$, counting
multiplicity, in $\{ z:|z-it/2|\leq t/2, |z| \geq 1\}$, and set
$$
\mathfrak{N} (r, f) = \int_1^r
\frac{ \mathfrak{n} (t, f) }{t^2} dt, \quad r \geq 1.
$$
The Tsuji characteristic is defined as
$$
\mathfrak{T} (r, f)  =
\mathfrak{m} (r, f)  +
\mathfrak{N} (r, f),
$$
where
$$\mathfrak{m} (r, f) =
\frac1{2 \pi} \int_{ \sin^{-1} (1/r)}^{ \pi - \sin^{-1} (1/r)}
\frac{ \log^+ | f(r \sin \theta e^{i \theta } )|}
{r \sin^2 \theta } d \theta.$$
The upper half-plane is thus exhausted by circles of diameter $r \geq 1$
tangent to the real axis at $0$. For non-constant $f$ and any
$a \in \C$ the first fundamental theorem then reads
\cite{GO,Tsuji0}
\begin{equation}
\mathfrak{T} (r, f) =
\mathfrak{T} (r, 1/(f - a)) + O(1), \quad r \to \infty,
\label{t1}
\end{equation}
and the lemma on the logarithmic derivative
\cite[p. 332]{LeO} gives
\begin{equation}
\mathfrak{m} (r, f'/f) = O(\log r + \log^+ \mathfrak{T} (r, f) )
\label{t2}
\end{equation}
as $r \to \infty $ outside a set of finite measure. Further,
$\mathfrak{T} (r, f)$ differs from a non-decreasing function
by a bounded additive term
\cite{Tsuji0}. Standard inequalities give
\begin{equation}
\mathfrak{T} (r, f_1 + f_2)
\leq \mathfrak{T} (r, f_1 ) +    \mathfrak{T} (r, f_2 ) +    \log 2,
\quad
\mathfrak{T} (r, f_1  f_2) \leq \mathfrak{T} (r, f_1 )
+    \mathfrak{T} (r, f_2 ),
\label{t3}
\end{equation}
whenever $f_1, f_2$ are meromorphic in $\overline{H}$.
Using the obvious fact that $\mathfrak{T}(r,1/z)=0$ for
$r\geq 1$
we easily derive from (\ref{t1})
and (\ref{t3})
that $\mathfrak{T}(r,f)$ is bounded if $f$ is a rational function.

A key role will be played by
the following two results from \cite{LeO}.
The first is obtained by a change of variables in
a double integral \cite[p. 332]{LeO}.

\begin{lem}\label{lemt1}
Let $Q(z)$ be meromorphic in $\overline{H}$, and for $r \geq 1$ set
\begin{equation}
m_{0\pi} (r, Q) =
\frac1{2\pi} \int_0^\pi \log^+ |Q(r e^{i \theta  } )| d \theta.
\label{t4}
\end{equation}
Then for $R \geq 1$ we have
\begin{equation}
\int_R^\infty \frac{m_{0\pi} (r, Q) }{r^3} dr \leq
\int_R^\infty \frac{\mathfrak{m} (r, Q) }{r^2} dr.
\label{t5}
\end{equation}
\end{lem}

The second result from \cite{LeO}
is the analogue for the half-plane of
Hayman's Theorem~3.5 from \cite{Haym1}.

\begin{lem}\label{lemt2}
Let $k \in \N$ and
let $f$ be meromorphic in $\overline{H}$, with
$f^{(k)} \not \equiv 1$. Then
$$
\mathfrak{T} (r, f)  \leq
\left(2 + \frac{1}{k}\right) \mathfrak{N} \left(r, \frac{1}{f}\right) +
\left(2 + \frac{2}{k}\right) \mathfrak{N}
\left(r, \frac{1}{f^{(k)} - 1}\right) +
O(\log r + \log^+\mathfrak{T} (r, f))
$$
as $r \to \infty $ outside a set of finite measure.
\end{lem}
Lemma \ref{lemt2} is established by following Hayman's proof exactly,
but using the Tsuji characteristic and the lemma on the logarithmic
derivative (\ref{t2}).

We also need the following result of Yong Xing Gu
(Yung-hsing Ku, \cite{Ku}). 
\begin{lem}\label{lemku} 
For every $k\in\N$, the meromorphic functions $g$ in an arbitrary domain
with the properties that $g(z)\neq 0$ and $g^{(k)}(z)\neq 1$
form a normal family. 
\end{lem}
A simplified proof of this result
is now available \cite{Zalc}.
It is based on a rescaling lemma
of Zalcman--Pang \cite{Pang} which permits an easy derivation of
Lemma \ref{lemku}
from the following result of Hayman:
{\em Let $k\in \N$ and let $g$ be a meromorphic function in the plane
such that
$g(z)\neq 0$ and $g^{(k)}(z)\neq 1$ for $z\in\C$. Then 
$g=\const$}, see \cite{Haym} or 
\cite[Corollary of Thm 3.5]{Haym1}.  

\section{Proof of Theorem \ref{thm1}}

Let $L,\psi,\phi$ be as in the hypotheses,
and assume that $\phi$ is transcendental but $L+L'/L$ has only
finitely many non-real zeros.
Condition (b) implies the Carath\'eodory inequality:
\begin{equation}
\frac{1}{5}|\psi(i)|\frac{\sin\theta}{r}<|\psi(re^{i\theta}
)|<
5|\psi(i)|\frac{r}{\sin\theta}, \quad r \geq 1,\quad
\theta\in (0,\pi),
\label{3}
\end{equation}
see, for example, \cite[Ch. I.6, Thm $8'$]{Le}.
\begin{lem}\label{lolem1}
The Tsuji characteristic of $L$ satisfies
$\mathfrak{T}(r,L) = O(\log r)$ as $r\to\infty$.
\end{lem}

\noindent{\em Proof.}
We apply Lemma \ref{lemt2} almost exactly as in
\cite[p. 334]{LeO}.  Let $g_1 = 1/L$. Then
$$
g_1' = - L' /L^2.
$$
Since $L$ has only real poles and since
$L + L'/L$ has by assumption finitely many non-real zeros it
follows that $g_1$ and $g_1' - 1$ have finitely many zeros in $H$.
Lemma \ref{lemt2} now gives
$\mathfrak{T} (r, g_1 ) = O( \log r )$ initially outside
a set of finite measure, and hence without exceptional set since
$\mathfrak{T} (r, g_1 ) $ differs from a non-decreasing function
by a bounded term. Now apply (\ref{t1}).
\hfill$\Box$

\begin{lem}\label{lolem}
The function
$\phi$ has order at most $1$.
\end{lem}

\noindent{\em Proof.}
Again, this
proof is almost identical to the corresponding argument in 
\cite{LeO}.
Lemmas \ref{lemt1}
and \ref{lolem1} give
$$
\int_R^\infty \frac{m_{0\pi}(r, L)}{r^3} dr
\leq  \int_R^\infty \frac{\mathfrak{m} (r, L)}{r^2} dr
= O( R^{-1} \log R ),
\quad R  \to \infty.
$$
Since
$ m_{0\pi}(r, 1/\psi ) = O( \log r ) $
by (\ref{3}), we obtain using (\ref{1})
$$
\int_R^\infty \frac{m_{0\pi}(r, \phi)}{r^3} dr = O( R^{-1} \log R ),
\quad R  \to \infty.
$$
But $\phi$ is entire and real on the real axis and so
$$
\| \log^+|\phi(re^{i\theta})|\| = 2 m_{0\pi} (r, \phi ).
$$
Since $\| \log^+|\phi(re^{i\theta})|\|$ is a non-decreasing
function of $r$ we deduce that
$$
\| \log^+|\phi(Re^{i\theta})|\|
= O( R \log R ), \quad R \to \infty,
$$
which proves the lemma.
\hfill$\Box$

\begin{lem}\label{lem2}
Let $\delta_1 > 0$ and $K > 1$. Then we have
\begin{equation}
\left| wL(w) \right| > K,
\quad |w| = r,
\quad
\delta_1 \leq \arg w \leq \pi -
\delta_1,
\label{9}
\end{equation}
for all $r$ outside a set $E_1$ of zero logarithmic density.
\end{lem}

\noindent{\em Proof.}
Choose $\delta_2$ with
$0 < \delta_2 < \delta_1$. Let
$$
\Omega_0 = \{ z \in \C : \frac12 < |z| < 2,
\quad
\frac{\delta_2}{2} < \arg z <
\pi - \frac{\delta_2}{2} \}.
$$
For $r \geq  r_0$, with
$r_0$ large, let $g_r(z) = 1/(rL(rz))$. Then $g_r(z) \neq 0$ on
$\Omega_0 $, since all poles of $L$ are real. Further,
$$
g_r'(z)= -L'(rz)/L(rz)^2.
$$
Since $L$ is analytic in $H$ and $L+L'/L$ has finitely many zeros in $H$
it follows that provided $r_0$ is
large enough the equation $g_r'(z) = 1$ has no solutions in
$\Omega_0$.
Thus the functions $g_r(z)$ form a
normal family on $\Omega_0$, by Lemma \ref{lemku}. 

Suppose that $|w_0| = r \geq  r_0$, and
$\delta_1 \leq \arg w_0 \leq
\pi - \delta_1$, and that
\begin{equation}\label{new}
\left| w_0 L(w_0) \right| \leq K.
\end{equation}
Then
$$
\left| g_r (z_0) \right| \geq 1/K, \quad z_0 =
\frac{w_0}{r},
$$
and so since the $g_r$ are zero-free and form a normal family
we have
\begin{equation}
\left| g_r (z) \right| \geq 1/K_1,
\quad |z| = 1, \quad
\delta_2 \leq \arg z \leq \pi -
\delta_2,
\label{11}
\end{equation}
for some positive constant
$K_1 = K_1 ( r_0, \delta_1, \delta_2, K)$, independent of $r$.  By
(\ref{1}),
(\ref{3}),
and
(\ref{11}) we have, for $|w| = r,\;
\delta_2 \leq \arg w \leq \pi -
\delta_2 $, the estimates
\begin{equation}
\left| wL(w) \right| =
| w \psi (w) \phi (w) |  \leq K_1, \quad
|  \phi (w) |  \leq K_2 = \frac{5 K_1}{| \psi (i) | \sin \delta_2 } .
\label{12}
\end{equation}
Thus (\ref{new}) implies (\ref{12}).
For $t \geq r_0$ let
$$
E_2(t) = \{ w \in \C : |w| = t,
| \phi (w) | > K_2 \}.
$$
Further, let $\theta (t)$ be the angular measure of $E_2(t)$,
and as in Lemma \ref{sub} let
$\theta^* (t) = \theta (t)$, except that
$\theta^* (t) = \infty $ if $E_2(t) = C(0, t)$.
Let
$$
E_3 = \{ t \in [r_0, \infty ) : \theta (t) \leq 4 \delta_2 \}.
$$
Since (\ref{new}) implies (\ref{12}), we have
(\ref{9}) for $t \in [r_0, \infty ) \setminus E_3$. Applying
Lemma~\ref{sub}
we obtain, since $\phi$ has order at most $1$ by Lemma \ref{lolem},
$$
(1 + o(1)) \log r \geq
\int_{r_0}^r \frac{ \pi dt }{ t \theta^* (t) } \geq
\int_{[{r_0}, r] \cap E_3 }  \frac{ \pi dt }{ 4 \delta_2 t  },
$$
from which it follows that $E_3$ has upper logarithmic density at
most $4 \delta_2 / \pi $. Since
$\delta_2$ may be chosen arbitrarily small, the lemma is
proved.
\hfill$\Box$
\vspace{.1in}

The estimates (\ref{3}) and
(\ref{9}) and the fact that $\phi$ is real now give
$$
| \phi (z) | >    \frac{K \sin \delta_1 }{5 | \psi (i) |  r^2}, \quad
\delta_1 \leq | \arg z | \leq
\pi - \delta_1,
$$
for $|z| = r$ in a set of logarithmic
density $1$. Since $\phi$ has order at most $1$ but is
transcendental, we deduce
(compare \cite[pp. 500-501]{HelW2}):

\begin{lem}\label{lem22}
The function $\phi$ has infinitely
many zeros.
\hfill$\Box$
\end{lem}

Let
\begin{equation}
F(z) = z - \frac{1}{L(z)}, \quad
F'(z) = 1 + \frac{L'(z)}{L(z)^2} .
\label{6}
\end{equation}
Since $L$ has only real poles and
$L+L'/L$ has finitely many non-real zeros we obtain at once:

\begin{lem}\label{F'crit}
The function $F$ has finitely many critical points over
$\C \setminus \R$, i.e. zeros $z$ of $F'$ with $F(z)$ non-real.
\hfill$\Box$
\end{lem}

\begin{lem}\label{lem1}
There exists
$\alpha \in H$ with the property that $F(z) \to \alpha $ as
$z \to \infty $ along a path $\gamma_\alpha $ in $H$.
\end{lem}
Lemma \ref{lem1} is a refinement of Theorem 4 of \cite{SS}, and will be proved
in \S\ref{lem1pf}.

Now set
\begin{equation}
g(z) = z^2 L(z) - z = \frac{zF(z)}{z - F(z) }, \quad
h(z) = \frac1{ F(z) - \alpha },
\label{hdef}
\end{equation}
in which $\alpha$ is as in Lemma \ref{lem1}. Then $g$ is analytic in $H
\cup \{ 0 \}$ and
(\ref{t3}),
(\ref{6}) and Lemma \ref{lolem1} give
$$
\mathfrak{T}(r, g) +
\mathfrak{T}(r, h) = O( \log r ),
\quad r \to \infty.
$$
Hence Lemma \ref{lemt1} leads to
\begin{equation}
\int_1^\infty \frac{ m_{0\pi} (r, g) }{r^3} dr +
\int_1^\infty \frac{ m_{0\pi} (r, h) }{r^3} dr < \infty,
\label{h1}
\end{equation}
in which
$m_{0\pi} (r, g) $ and
$m_{0\pi} (r, h) $ are as defined in (\ref{t4}).

\begin{lem}\label{lem11}
The function $F$ has at most four finite non-real asymptotic values.
\end{lem}

\noindent{\em Proof.}
Assume the contrary. Since $F(z)$ is real on the real axis
we may take distinct finite non-real
$\alpha_0, \ldots, \alpha_n,\; n \geq 2,$
such that $F(z) \to \alpha_j$ as $z\to \infty$ along a simple path
$\gamma_j: [0, \infty ) \to H \cup \{0 \}$. Here we assume
that $\gamma_{j} (0) = 0$, that
$\gamma_{j} (t) \in H$ for $t > 0$, and that
$\gamma_{j} (t) \to \infty $ as $t \to \infty $. We may further
assume that
$\gamma_{j} (t) \neq
\gamma_{j'} (t')$ for $t > 0,\, t'>0,\, j \neq j'$.

Re-labelling if necessary, we obtain $n$ pairwise disjoint simply
connected domains $D_1, \ldots, D_n$ in $H$, with $D_j$ bounded
by
$\gamma_{j-1}$ and
$\gamma_j$, and for $t > 0$ we let
$\theta_j(t)$ be the angular measure of the intersection of
$D_j$ with the circle
$C(0, t)$. By (\ref{hdef}), the function $g(z)$
tends to $\alpha_j$ as $z \to \infty $ on $\gamma_j$, and
so $g(z)$ is unbounded on each $D_j$ but bounded on the finite boundary
$\partial D_j$ of each $D_j$. Let
$c $ be large and positive, and for each $j$ define
\begin{equation}
u_j(z) = \log^+ | g(z)/c |, \quad z \in D_j.
\label{ujdef}
\end{equation}
Set $u_j(z) = 0$ for $z \not \in D_j$. Then $u_j$ is continuous,
and subharmonic in the plane since $g$ is analytic in $H \cup \{ 0 \}$.

Lemma \ref{sub} 
gives, for some $R > 0$ and for each~$j$,
$$
\int_R^{r} \frac{\pi dt}{t \theta_j(t)} \leq
\log \|u_j(4re^{i\theta})\| + O(1)
$$
as $r \to \infty $. Since $u_j$ vanishes outside $D_j$ we deduce
using (\ref{ujdef}) that
\begin{equation}
\int_R^{r} \frac{\pi dt}{t \theta_j(t)}
\leq \log m_{0\pi}(4r, g ) + O(1), \quad r \to \infty,
\label{7}
\end{equation}
for all $j \in \{ 1, \ldots, n \}$. However, the Cauchy-Schwarz
inequality gives
$$
n^2 \leq \sum_{j=1}^n \theta_j(t) \sum_{j=1}^n \frac1{\theta_j(t) }
\leq \sum_{j=1}^n \frac{\pi}{\theta_j(t) }
$$
which on combination with (\ref{7}) leads to, for some positive constant
$c_3$,
$$
n \log r \leq \log m_{0\pi} (4r, g) + O(1), \quad
m_{0\pi} (r, g) \geq c_3 r^n, \quad r \to \infty.
$$
Since $n \geq 2$ this contradicts (\ref{h1}), and
Lemma \ref{lem11} is proved.
\hfill$\Box$
\vspace{.1in}

{}From Lemmas \ref{F'crit} and \ref{lem11}
we deduce
that the inverse function $F^{-1}$ has finitely many non-real singular
values. Using
Lemma \ref{lem1},
take $\alpha \in H$ such that $F(z) \to \alpha $ along a
path $\gamma_{\alpha}$ tending to infinity in $H$,
and take
$\varepsilon_0$ with
$0 < \varepsilon_0 < {\rm Im} (\alpha)$ such that $F$
has no critical or asymptotic values in
$0 < | w - \alpha | \leq \varepsilon_0$. Take a
component $C_0$
of the set $ \{ z \in \C : | F(z) - \alpha | < \varepsilon_0 \}$
containing an unbounded subpath of $\gamma_{\alpha}$. Then by
a standard argument \cite[XI.1.242]{Ne}
involving a logarithmic change of variables
the inverse function $F^{-1}$ has a logarithmic singularity over
$\alpha$,
the component $C_0$ is
simply connected, and
$F(z) \neq \alpha $ on $C_0$. Further, the boundary of $C_0$
consists of a single simple curve going to infinity in both
directions. Thus we may define a
continuous, non-negative, non-constant subharmonic function in the plane
by
\begin{equation}
u(z)=\log\left|\frac{\varepsilon_0}{F(z)-\alpha}\right|=
\log|\varepsilon_0h(z)|
\quad (z \in C_0),
\quad \quad u(z) = 0 \quad  (z \not \in C_0),
\label{8}
\end{equation}
using (\ref{hdef}).

The next lemma follows from (\ref{6})
and (\ref{8}).

\begin{lem}\label{lem3}
For large $z$ with $|z L(z)| > 3$ we have
$|F(z) - \alpha | > |z|/2$ and $u(z) = 0$.
\hfill$\Box$
\end{lem}

\begin{lem}\label{lem4}
We have
\begin{equation}
\lim_{r \to \infty }
\frac{ \log \| u(re^{i\theta})\| }
{ \log r } = \infty.
\label{13}
\end{equation}
\end{lem}

\noindent{\em Proof.}
Apply Lemma \ref{lem2}, with $K = 3$ and $\delta_1$ small and positive.
By Lemma~\ref{lem3} we have
$ u(z) = 0$ if
$\delta_1 \leq |\arg z | \leq \pi -
\delta_1 $ and
$|z|$ is large but not in $E_1$. For large $t$ let
$\sigma (t)$ be the angular measure of that subset of
$C(0, t)$ on which $u(z) > 0$. Since $u$ vanishes on the real axis
Lemma \ref{sub} and Lemma \ref{lem2} give, for some $R > 0$,
$$
\log \| u(4re^{i\theta}) \|+ O(1) \geq
\int_R^r \frac{ \pi dt}{ t \sigma (t) } \geq
\int_{ [R, r] \setminus E_1 }  \frac{ \pi dt}{  4 \delta_1 t  } \geq
\frac{\pi}{4 \delta_1} ( 1 - o(1) ) \log r
$$
as $r \to \infty $. Since $\delta_1$ may be chosen  arbitrarily small
the lemma follows.
\hfill$\Box$
\vspace{.1in}

Now (\ref{8}) and the fact that $u$ vanishes outside $C_0$ give
$$
 \| u(re^{i\theta})\| \leq
m_{0\pi} (r, h) + O(1),
$$
from which we deduce using (\ref{13}) that
$$
\lim_{r \to \infty}
\frac{ \log m_{0\pi} (r, h) }{ \log r } = \infty.
$$
This obviously contradicts (\ref{h1}), and Theorem \ref{thm1}
is proved.
\hfill$\Box$

\section{Proof of Lemma \ref{lem1}}\label{lem1pf}

The proof is based essentially on Lemmas 1 and 5 and Theorem 4 of
\cite{SS}.
Assume that there is
no $\alpha \in H$ such that $F(z)$ tends to
$\alpha $ along a path tending to infinity in $H$.

Let
$$
W = \{ z \in H : F(z) \in H \}, \quad
Y = \{ z \in H : L(z) \in H \}.
$$
Then $Y \subseteq W$, by (\ref{6}), so that each
component $C$ of $Y$ is contained in
a component $A$ of $W$.

\begin{lem}\label{lemSS1}
To each component $A$ of $W$ corresponds a finite number
$v(A)$ such that $F$ takes every value
at most $v(A)$ times in $A$ and has at most $v(A)$
distinct poles on $\partial A$.
\end{lem}

\noindent
{\em Proof.} Using Lemma \ref{F'crit}, cut $H$ along a simple
polygonal curve starting from $0$ so that the resulting simply connected
domain $D$ contains no critical value of $F$. Let
$X=\{ z\in H: F(z)\in D\}$.
By analytic continuation of the inverse function, every
component $B$ of
$X$ is simply connected and conformally
equivalent under $F$ to $D$.
Moreover, if the finite boundary
$\partial B$ contains no critical
point $z$ of $F$ with $F(z) \in H$
then the branch of the inverse
function mapping $D$ onto $B$
may be analytically continued throughout
$H$, and in this case
$F$ is univalent in the component $A$ of $W$ containing
$B$. 

Fix a component $A$ of $W$ on which $F$ is not univalent.
If $z\in A$ with $F(z)\in D$ then
$z$ lies in a component $B \subseteq A$ of $X$, such that $\partial B$
contains a critical point of $F$.
Since $F$ has finitely
many critical points over $H$, and every critical point
can belong to the boundaries of at most finitely many components
$B$ of $X$, it follows that $A$ contains
finitely many such components $B$.
Application of the open mapping theorem gives us a finite $v(A)$
satisfying the first statement of the lemma.

If $z_0\in \partial A$ is a pole
of $F$, then for an arbitrarily small
neighbourhood $N$ of $z_0$, $F$ assumes in $N\cap A$
all sufficiently large values in $H$. It follows that
there are no more than $v(A)$ distinct poles of $F$ on
$\partial A$.
\hfill$\Box$

\begin{lem}\label{lemSS2}
There are infinitely many components $A$ of $W$ such that
$A$
contains an unbounded component $C$ of $Y$
and $\partial A\cap\partial C$ contains a
zero of $L$. 
\end{lem}

\noindent{\em Proof.}
We first note \cite[Lemma 1]{SS} that $L$ has no poles in
the closure of $Y$.
To see this, let $x_0$ be a pole of $L$. Then $x_0$ is real
and is a simple pole
of $L$ with positive residue. Hence
$\lim_{y \to 0+} {\rm Im} (L(x_0 + iy)) = - \infty $
and since $L$ is univalent on an open disc
$N_0 = B(x_0, R_0)$ it follows that ${\rm Im}( L(z)) < 0$ on
$N_0 \cap H$.

Thus every component $C$ of $Y$ is unbounded by
the maximum principle, since ${\rm Im} (L(z))$ is harmonic in $H$
and vanishes on $\partial C$.

Next, we recall from Lemma \ref{lem22}
that $\phi$ has infinitely many zeros;
by the hypotheses (\ref{1})
and (c) these
must be zeros of $L$.
Suppose first that\\
\\
(I) $L$ has infinitely many non-real zeros.\\
\\
Then $L$ has infinitely many zeros $\eta \in H$, and
$F( \eta ) = \infty $. Each such $\eta$ lies on
the boundary of a component $C$ of $Y$, and $C$ is contained
in a component $A$ of $W$, and $\eta \in \partial A$. By Lemma
\ref{lemSS1} we obtain in this way infinitely many components
$A$.

Since $L$ is real on the real axis
we obtain the same conclusions if either of the following conditions
hold: \\
\\
(II) $L$ has infinitely many multiple zeros;\\
\\
(III) $L$ has infinitely many real zeros $x$ with $L'(x) > 0$.\\
\\
We assume henceforth that neither (I) nor (II) holds,
and will deduce (III).
Let $\{ a_k \}$ denote the poles of $L$, in increasing order.
Then there are two possibilities.
The first is that there exist infinitely many intervals
$(a_k, a_{k+1})$ each containing at least one zero $x_k$ of $\phi$.
Since $\psi$ must have negative residues by (b) there must
be a zero $y_k$ of $\psi$ in $(a_k, a_{k+1})$,
and we may assume that $y_k \neq x_k$, since $L$ has by assumption finitely
many multiple zeros.
But then the graph of $L$ must cut the real
axis at least twice in $(a_k, a_{k+1})$, and so there exists
a zero $x$ of $L$ in $(a_k, a_{k+1})$ with $L'(x) > 0$. Thus we obtain
(III).

The second possibility is that we have infinitely many pairs
of zeros $a, b$ of $\phi$ such that $L$ has no poles on
$[a, b]$. In this case we again obtain a zero $x$ of $L$
with $L'(x) > 0$,
this time in $[a, b]$, and again we have (III).
\hfill$\Box$
\vspace{.1in}

We now complete the proof of Lemma \ref{lem1}.
Combining Lemmas \ref{F'crit},
\ref{lemSS1} and
\ref{lemSS2} we obtain at least one zero $\eta $ of $L$,
with $\eta \in \partial A \cap \partial C$, in which
$A, C$ are components of $W, Y$ respectively,
$C$ is unbounded and $C \subseteq A$,
and $F$ has no critical point in $A$. It then follows
by analytic continuation of the inverse function in $H$ that
$F$ maps $A$ univalently onto $H$. Since $F$ takes near $\eta$ all values
$w$ of positive imaginary part and
large modulus, it follows that
$F(z)$ is bounded as $z \to \infty $ in $A$, so that
$L(z) \to 0$ as $z \to \infty $ in $A$, and hence
as $z \to \infty $ in $C$. This contradicts the maximum principle.
\hfill$\Box$

\vspace{.1in}

{\em W.B.: 
Mathematisches Seminar,

Christian--Albrechts--Universit\"at zu Kiel,

Ludewig--Meyn--Str.\ 4,

D--24098 Kiel,

Germany

{bergweiler@math.uni-kiel.de}

\vspace{.1in}

A.E.: Purdue University, 

West Lafayette IN 47907

USA

eremenko@math.purdue.edu
\vspace{.1in}

J.K.L.: School of Mathematical Sciences,

University of Nottingham,

Nottingham NG7 2RD

UK

jkl@maths.nottingham.ac.uk}
\end{document}